\documentclass[a4paper]{amsart}
\usepackage[utf8]{inputenc}
\usepackage[english]{babel}

\usepackage{amscd,amssymb}
\usepackage{amsmath}
\usepackage{amsthm}
\usepackage{graphicx}
\usepackage{fancybox}
\usepackage{epic,eepic}
\usepackage{amstext}
\usepackage{mathtools}
\usepackage{esint}
\usepackage[square,numbers]{natbib}
\usepackage{booktabs}
\usepackage[hide links]{hyperref}
\usepackage{comment}
\usepackage{mathtools}
\usepackage{tikz,pgfplots}

\usepackage[noabbrev, capitalise, nameinlink]{cleveref}
\crefname{equation}{}{}
\crefname{assumption}{Assumption}{Assumptions}
\crefformat{equation}{\textup{#2(#1)#3}}

\newtheorem{theorem}{Theorem}[section]
\newtheorem{corollary}[theorem]{Corollary}
\newtheorem{lemma}[theorem]{Lemma}

\theoremstyle{definition}

\theoremstyle{remark}
\newtheorem{remark}[theorem]{Remark}
\numberwithin{theorem}{section}
\numberwithin{equation}{section}
\numberwithin{figure}{section}

\def\Th{\mathcal{T}_h}

\def\diam{\operatorname{diam}}

\def\with{\,:\,}

\DeclareMathOperator*{\argmin}{arg\,min}

\newcommand{\tnof}[1]{\| #1 \|_{L^2}}

\newcommand{\tspof}[2]{{( #1,#2 )}_{L^2}}

\numberwithin{equation}{section}
\numberwithin{theorem}{section}

\AtBeginDocument{%
	\def\MR#1{}
}

\begin{document}
	
\title[Positivity preserving FEM for the Gross-Pitaevskii ground state]{Positivity preserving finite element method for the Gross-Pitaevskii ground state: discrete uniqueness and global convergence}
\author[M.~Hauck, Y.~Liang, D.~Peterseim]{Moritz Hauck$^*$, Yizhou Liang$^\dagger$, Daniel Peterseim$^{\dagger,\ddagger}$}
\address{${}^+$ Department of Mathematical Sciences, Chalmers University of Technology and University of Gothenburg, 41296 Göteborg, Sweden}
\email{hauck@chalmers.se}
\address{${}^{\dagger}$ Institute of Mathematics, University of Augsburg, Universit\"atsstr.~12a, 86159 Augsburg, Germany}
\address{${}^{\ddagger}$Centre for Advanced Analytics and Predictive Sciences (CAAPS), University of Augsburg, Universit\"atsstr.~12a, 86159 Augsburg, Germany}
\email{\{yizhou.liang, daniel.peterseim\}@uni-a.de}
\thanks{The work of M.~Hauck is supported by the Knut and Alice Wallenberg foundation postdoctoral program in mathematics for researchers from outside Sweden (Grant No. KAW 2022.0260). The work of Y.~Liang is funded by a Humboldt Research Fellowship for postdocs from the Alexander von Humboldt Foundation. D.~Peterseim is supported by the European Research Council ERC under the European Union's Horizon 2020 research and innovation program (RandomMultiScales, Grant Agreement No.~865751). 
}
	
\begin{abstract}
We propose a positivity preserving finite element discretization for the nonlinear Gross-Pitaevskii eigenvalue problem. The method employs mass lumping techniques, which allow to transfer the uniqueness up to sign and positivity properties of the continuous ground state to the discrete setting. We further prove that every non-negative discrete excited state up to sign coincides with the discrete ground state. This allows one to identify the limit of fully discretized gradient flows, which are typically used to compute the discrete ground state, and thereby establish their global convergence. Furthermore, we perform a rigorous a priori error analysis of the proposed non-standard finite element discretization, showing optimal orders of convergence for all unknowns. Numerical experiments illustrate the theoretical results of this paper.
\end{abstract}



\maketitle
	
\section{Introduction}
The Gross-Pitaevskii eigenvalue problem arises in quantum physics where it describes the stationary quantum states of Bose-Einstein condensates. The problem involves a non-negative confinement (or trapping) potential $V \in L^\infty(\Omega)$ and a positive parameter $\kappa$ describing the repulsive interaction of the particles in the condensate. As domain we consider a bounded convex Lipschitz domain ${\Omega \subset \mathbb R^d}$, ${d \in \{1,2,3\}}$, noting that the restriction to a (sufficiently large) domain along with homogeneous Dirichlet boundary conditions is a reasonable modeling assumption for quantum states at the lowest part of the energy spectrum. The Gross-Pitaevskii eigenvalue problem then seeks  $L^2$-normalized eigenstates $\{u_j \with j \in \mathbb N\}\subset H^1_0(\Omega)$ and corresponding eigenvalues $\lambda_j \in \mathbb R$ such that
\begin{equation}\label{eq:evpstrong}
	-\Delta u_j +Vu_j +\kappa |u_j|^2u_j = \lambda_j u_j
\end{equation}
holds in the weak sense. The function $|u_j|^2$ represents the density of the stationary quantum state $u_j$ and $\lambda_j$ denotes the corresponding chemical potential. All eigenvalues of~\cref{eq:evpstrong} are real and positive and the smallest eigenvalue is simple, cf.~\cite{CCM10}. In the following, we assume without loss of generality that the ordering of the eigenvalues is non-decreasing, i.e., we have that $0<\lambda_1< \lambda_2\leq \dots\;$.

 Eigenvalue problem \cref{eq:evpstrong} can be considered as the Euler-Lagrange equation for critical points of the Gross-Pitaevskii energy, defined for all $v \in H^1_0(\Omega)$~as
\begin{equation}
	\label{eq:energy}
			\mathcal{E}(v) \coloneqq \tfrac12\tspof{\nabla v}{\nabla v} + \tfrac12 \tspof{Vv}{v} + \tfrac{\kappa}{4} \tspof{|v|^2 v}{v}.
\end{equation}
The ground state, which is the stationary quantum state of lowest energy, can be characterized as the solution to the following constrained minimization problem
	\begin{equation}
	\label{eq:gs}
	u \in  \argmin_{v \in H^1_0(\Omega)\with \tnof{v} = 1} \mathcal E (v).
\end{equation}
There are several important theoretical results that hold for the ground state.  First, the ground state exists and is unique up to sign. Second, the ground state coincides, up to the sign, with the eigenfunction~$u_1$ of \cref{eq:evpstrong} corresponding to the smallest eigenvalue $\lambda_1$. Third, the ground state satisfies $|u|>0$ in $\Omega$, which means that the sign of the ground state can be chosen such that it is positive in~$\Omega$. Moreover, by Picone's inequality~\cite{Brasco2014}, any eigenfunction~$u_j$ without a change of sign must necessarily coincide with the ground state. The proofs of these theoretical results can be found, e.g., in \cite{CCM10,HenP20}.

In the literature there are a number of spatial discretizations to approximate the Gross-Pitaevskii ground state. Such discretizations may be based, for example, on finite element methods (FEMs) in primal \cite{Zho04,CCM10,CHZ11} or mixed \cite{GHL24} formulations, spectral and pseudospectral methods~\cite{CCM10,Bao13}, or mesh-adaptive methods~\cite{Dan10,He21}. Multiscale methods such as~\cite{HMP14b,MR4379968,henning2023optimal,PWZ23}, which are based on the (Super-)Localized Orthogonal Decomposition (cf.~\cite{MaP14,HeP13,HaPe21b}), are also popular discretization methods. 
We emphasize that very little is known about whether the discrete solutions obtained by these methods satisfy the above-mentioned properties of the continuous ground state. Very recently, partial progress has been made in \cite{positivitygpe}, where the  uniqueness and positivity of the discrete ground state for a lumped finite difference discretization was proved. To the best of our knowledge, it is still an open question whether the other properties can be transferred to the discrete setting. 

In addition to the choice of a suitable spatial discretization, one has to deal with the solution of the resulting  finite-dimensional constrained minimization problem, which is a discrete version of~\cref{eq:gs}. 
An overview of algorithms for this purpose is given in the recent review paper~\cite{Henning2024}. 
Popular methods include Sobolev gradient flows  such as~\cite{BaD04,RSS09,DaK10,KaE10}, Riemannian optimization methods \cite{DanP17,ALT17,AltmannPeterseimStykel2022,altmann2023riemannian}, Newton-type algorithms such as the (shifted) $J$-method \cite{JarKM14,AltHP21}, and the self-consistent field iteration~\cite{Can00,DioC07}. We highlight~\cite{HenP20}, where a gradient flow method based on an energy-adaptive metric is proposed. After appropriate pseudo-time discretization, the resulting iteration was reinterpreted as an energy-adaptive Riemannian gradient descent method in~\cite{AltmannPeterseimStykel2022}. The global convergence of the iteration to the continuous ground state was shown in~\cite{HenP20}. Crucial ingredients to prove this convergence result are the uniqueness and positivity of the ground state as well as Picone's inequality to identify the limit of the iteration as the ground state, given a non-negative initial guess. These properties no longer hold after discretization in space using, e.g., standard  methods such as linear or quadratic~FEMs. Therefore, the reasoning from the continuous setting cannot be transferred to the discrete setting, and a global convergence result for the fully discrete case is still unknown.

To address the lack of positivity preservation, inspired by work on discrete maximum principles for convection-diffusion problems (see, e.g., the recent review article \cite{Barrenechea2024}), we consider a  mass-lumped linear FEM for the Gross-Pitaevskii eigenvalue problem. For this particular discretization, we are able to prove that it preserves the above-mentioned properties of the continuous ground state. More precisely, we show that the discrete ground state is unique up to sign and coincides with the discrete eigenfunction corresponding to the smallest discrete eigenvalue. We also prove that the discrete ground state is positive in $\Omega$ and that a discrete version of Picone's inequality holds. These results allow us to establish the global convergence of a fully discretized version of the Sobolev gradient flow of \cite{HenP20} to the discrete ground state, for any non-negative initial guess.  
 
Furthermore, a rigorous a priori analysis of the proposed non-standard FEM is performed, proving optimal orders of convergence as the mesh size is decreased. Note that the achieved orders are the same as for the standard linear FEM, which does not preserve positivity.    We emphasize that for estimating the lumping errors, the corresponding theory for linear eigenvalue problems cannot be directly applied, and more specific tools need to be developed. 
The theoretical results of this paper are supported by a number of numerical experiments. The corresponding code is available at \url{https://github.com/moimmahauck/GPE_P1_lumped}. 

The paper is organized as follows: In \cref{sec:discretization} we introduce the proposed discretization for the Gross-Pitaevskii problem. The uniqueness and positivity of the discrete ground state are proved in \cref{sec:uniquenesspositivity}. A discrete version of Picone's inequality is derived in \cref{sec:globconv}, which is then used to establish the global convergence of a fully discretized Sobolev gradient flow to the discrete ground state. An a priori error analysis of the proposed method is given in \cref{sec:erroranalysis}. Finally, in \cref{sec:numexp} we provide numerical experiments to support our theoretical findings.
	
\section{Finite element discretization}\label{sec:discretization}

Consider a geometrically conforming and shape-regular hierarchy of simplicial finite element meshes $\{\Th\}_{h}$ of the domain $\Omega$. We denote the elements of a mesh~$\Th$ in the hierarchy by $K$ and define the mesh size~$h>0$ as the maximum diameter of the elements in~$\Th$, i.e., $h \coloneqq \max_{K \in \Th} \diam(K)$. Given the mesh $\mathcal T_h$, we denote by $n \in \mathbb N$ its number of interior and boundary nodes and by $\{p_j\with j = 1,\dots,n\}$ the coordinates of the nodes. For the discretization of the Gross-Pitaevskii problem, we use a linear FEM combined with a classical mass lumping approach, cf.~\cite{XuZ99}. Henceforth, we denote by $V_h$ the ansatz space of the linear FEM consisting of globally continuous $\mathcal T_h$-piecewise polynomials of total degree at most one, and by~$V_h^0$ the subspace satisfying homogeneous Dirichlet boundary conditions at $\partial \Omega$. We define the mass-lumped bilinear form for $\mathcal T_h$-piecewise continuous functions $v,w$ as
\begin{equation}
	\label{eq:lumpedbilinearform}
	\ell(v,w) \coloneqq \sum_{K\in\mathcal{T}_h} \ell_K(v,w),\qquad \ell_K(v,w)\coloneqq \frac{|K|}{d+1}\sum_{j=1}^{d+1}v|_{K}(p_{\tau_K(j)})w|_{K}(p_{\tau_K(j)}),
\end{equation}
where $|K|$ denotes the volume of the simplex $K$ and $\tau_K\colon \mbox{\{1,\dots,d+1\}}\to \{1,\dots,n\}$ maps the local node indices of the element $K$ to the corresponding global node indices. For discrete functions $v_h \in V_h$, this bilinear form is actually an inner product, and we denote its induced norm by 
\begin{equation*}
	\|v_h\|_{\ell}^2: = \ell(v_h,v_h).
\end{equation*}

The proposed method is based on an energy functional obtained by replacing the $L^2$-inner products in the definition of the energy $\mathcal E$, cf.~\cref{eq:energy}, by their lumped counterparts. Assuming that the potential $V$ is $\mathcal T_h$-piecewise continuous, the resulting energy functional is for all $v_h \in V_h^0$ defined~as 
\begin{equation}
	\label{eq:discenergy}
	\mathcal E_h(v_h) \coloneqq \tfrac12 \tspof{\nabla v_h}{\nabla v_h} + \tfrac12 \ell(Vv_h,v_h) + \tfrac{\kappa}{4} \ell(|v_h|^2 v_h,v_h).
\end{equation}
The discrete ground state $u_h \in V_h^0$ of the proposed method is then defined as the solution to the finite-dimensional constrained minimization problem
\begin{equation}
	\label{eq:gsdisc}
	u_h \in \argmin_{v_h \in V_h^0 \with \|v_h\|_\ell = 1} \mathcal E_h(v_h).
\end{equation}
In the discrete setting, the boundedness of the norms of the minimizing sequence implies the strong convergence of a subsequence (Bolzano-Weierstrass theorem). Thus there always exist the discrete energy minimizers $u_h$ and $-u_h$. Note that the uniqueness up to sign of the minimizer, which holds for the continuous problem, is generally not clear in the discrete setting. The approach to prove uniqueness in the continuous setting is to transform the original non-convex problem into a convex one acting on the densities (the squared modulus of the state). This approach is not applicable for standard FEMs, because the set 
\begin{equation*}
	\{v_h^2\with v_h \in V_h^0,\; \|v_h\|_{L^2} = 1\}
\end{equation*} 
is not convex, cf.~\cite{CCM10}. For non-standard FEMs, such as non-conforming or mixed FEMs, the latter lack of convexity can be overcome in some cases. However, the lack of certain positivity properties still remains a problem. Therefore, we adopt a mass lumping approach similar to the one used in the context of positivity preservation for convection-diffusion problems; see also the review article \cite{Barrenechea2024}. 

The Euler-Lagrange equations for  critical points of the energy $\mathcal E_h$ give rise to the  nonlinear eigenvalue problem: Seek $(v_h,\mu_h)\in V_h^0 \times \mathbb{R}$ with $\|v_h\|_\ell =1$ such that  
\begin{equation}\label{eq:discrete_eigen}
	(\nabla v_h,\nabla w_h)_{L^2} + \ell(Vv_h,w_h) + \kappa \ell(|v_h|^2v_h,w_h) = \mu_h\ell(v_h,w_h)
\end{equation}
holds for all $w_h\in V_h^0$. 
Note that the discrete ground state $u_h$ is an eigenfunction of~\cref{eq:discrete_eigen}, and we denote the corresponding eigenvalue by $\lambda_h$. We emphasize that, in contrast to the continuous setting, it is generally not clear that the discrete ground state coincides up to sign  with the eigenfunction corresponding to the smallest eigenvalue, and that~$\lambda_h$ is a simple eigenvalue, cf.~\cite{CCM10}. 

\section{Uniqueness and positivity of the discrete ground state}\label{sec:uniquenesspositivity}

In this section, we will show the uniqueness and positivity of the discrete ground state obtained by the proposed mass-lumped FEM. In addition, we will prove that the ground state eigenvalue is the smallest eigenvalue of  eigenvalue problem~\cref{eq:discrete_eigen}~and that it is simple. These results are not only of physical interest, but also lay the foundation for the proof of a discrete Picone-type inequality in \cref{sec:globconv}. This inequality is essential for establishing the global convergence of fully discretized Sobolev gradient flows to the discrete ground state.

To derive the desired discrete uniqueness and positivity properties we need to impose certain geometric conditions on the mesh $\mathcal T_h$. More precisely, denoting by ${\{\Lambda_j\with j = 1,\dots,m\}\subset V_h^0}$ the set of hat functions corresponding to the interior nodes of the mesh $\mathcal T_h$, where~$m \in \mathbb N$ is the number of interior nodes, one needs to ensure that the stiffness matrix $\mathbf S \in \mathbb R^{m\times m}$ with $S_{ij}\coloneqq (\nabla \Lambda_j,\nabla \Lambda_i)_{L^2}$ is an M-matrix, cf.~\cite{Plemmons1977}. The M-matrix property is classical in the context of discrete maximum principles, and various sufficient geometric conditions on the mesh $\mathcal T_h$ have been identified in the literature. In two dimensions, the M-matrix property is satisfied under the condition that the sum of the angles opposite to any edge are less than or equal to $\pi$, which is closely related to $\mathcal T_h$ being a Delaunay triangulation, cf.~\cite{XuZ99}. For the three-dimensional case, more restrictive conditions are typically imposed. One may, e.g., consider non-obtuse tetrahedral meshes for which the M-matrix property is known to hold, cf.~\cite{KrP00}.  For a discussion of appropriate refinement strategies for tetrahedra, we refer to~\cite{Brandts2020} and the references therein. In addition to the M-matrix property, we will make the technical assumption that~$\mathbf S$ is irreducible. This assumption is typically not restrictive and, if not already satisfied, can be ensured by appropriate local refinement of the considered~mesh.

The following theorem encapsulates the first major result of this paper. By deriving a strictly convex minimization problem for $|u_h|^2$, we are able to prove the desired  uniqueness and positivity properties of the discrete ground state. 
We emphasize that, in contrast to \cite{positivitygpe}, where a similar result is proved for a mass-lumped finite difference discretization, our proof does not rely on the explicit knowledge of the stiffness matrix and thus also allows the consideration of unstructured meshes. Furthermore, also the techniques used in the proofs are different. Contrary to~\cite{positivitygpe}, our proof does not make explicit use of the Perron-Frobenius~theorem.
 
\begin{theorem}[Uniqueness and positivity of discrete ground state]\label{theo:uni_posi:mini}
	Suppose that the stiffness matrix $\mathbf S$ is an irreducible M-matrix. Then the discrete ground state~$u_h$ defined in~\cref{eq:gsdisc} is unique up to sign. Furthermore, by appropriately flipping its sign, the discrete ground state can be chosen to be positive in~$\Omega$.
\end{theorem}

\begin{proof}
	This proof is done in two steps: In Step 1 we prove the positivity of the discrete ground state and in Step 2 its uniqueness. Below, the proof of~\cite[Thm.~2.1]{Lieb2000} will serve as source of inspiration. There the uniqueness and positivity of the continuous ground state was first proved; see also \cite{CCM10}.
	
	\emph{Step 1:} 
	The discrete ground state $u_h \in V_h^0$ can be written as the linear combination $u_h = \sum_{j=1}^m u_j\Lambda_j$, where $\mathbf u \coloneqq (u_j)_{j = 1}^m \in \mathbb R^m$ denotes the corresponding coefficient vector. Denoting a non-negative version of the discrete ground state by ${\bar{u}_h = \sum_{j=1}^m|u_j|\Lambda_j}$, it can be shown that $\mathcal{E}_h(\bar{u}_h)\leq \mathcal{E}_h(u_h)$. The proof of this inequality exploits  that the off-diagonal entries of $\mathbf S$ are non-positive, which holds by the M-matrix property of $\mathbf S$. Therefore, we can assume without loss of generality that~$u_h \geq 0$. 
	Next, we will show that $u_h>0$ in $\Omega$. Given a $\mathcal T_h$-piecewise continuous non-negative weighting function $w$, we denote by $\mathbf w = (w_j)_{j = 1}^k\in \mathbb R^{k}$ with $k\coloneqq (d+1)\cdot \#\mathcal T_h$ the corresponding vector of element-wise nodal evaluations. We define the weighted lumped mass matrix $\mathbf M(\mathbf w)\in \mathbb R^{m\times m} $ by $M(\mathbf w)_{ij}\coloneqq \ell(w\varphi_j,\varphi_i)$ and write~$\mathbf M$ for the unweighted lumped mass matrix. Note that for any $\mathbf w \in \mathbb R^k$, the resulting weighted lumped mass matrix $\mathbf M(\mathbf w)$ is diagonal. Furthermore, we denote by $\mathbf P \in \mathbb R^{k\times m}$ the canonical prolongation matrix. Denoting by $\mathbf V \in \mathbb R^k$ the vector of element-wise nodal evaluations of $V$, and by $\mathbf u^2$ the component-wise square of $ \mathbf u$, the coefficient vector $\mathbf u$ solves the following generalized eigenvalue problem
	\begin{equation}
		\label{eq:discevp}
		\mathbf A\mathbf u = \lambda_h \mathbf M \mathbf u,\qquad \mathbf A \coloneqq \mathbf S + \mathbf M(\mathbf V) + \mathbf M(\mathbf P\mathbf u^2),
	\end{equation}
	cf.~\cref{eq:discrete_eigen}. This is a linear problem because we fixed the vector $\mathbf u$ in the definition of the matrix $\mathbf A$.
	Since the matrix $\mathbf A$ is the sum of an M-matrix and a non-negative diagonal matrix, it is also an M-matrix, cf.~\cite[Thm.~2]{Plemmons1977}. Additionally, the matrix~$\mathbf A$ is irreducible. This is because $\mathbf S$ is irreducible and has positive diagonal entries, and therefore adding a non-negative diagonal matrix does not change the matrix' sparsity pattern. Introducing the variable $\mathbf v \coloneqq \mathbf M^{-1/2}\mathbf u$ allows us to write~\cref{eq:discevp} as the classical eigenvalue problem $\mathbf B \mathbf v = \lambda_h \mathbf v$ with $\mathbf B \coloneqq \mathbf M^{1/2}\mathbf A \mathbf M^{1/2}$. Since $\mathbf M$ is a diagonal matrix with positive diagonal entries, the matrix $\mathbf B$ is symmetric positive definite. A consequence of the M-matrix property of $\mathbf A$ is that $\mathbf A^{-1} \geq 0$ holds entry-wise. It is easy to verify that also  $\mathbf B^{-1} \geq 0$ holds entry-wise. Due to the irreducibility of $\mathbf B$ and hence also its inverse, there exists $k \in \mathbb N$ such that $\mathbf B^{-k}>0$ holds entry-wise. Since $\mathbf B^{-k}\mathbf v = \lambda_h^{-k}\mathbf v$, we can conclude that $\mathbf v >0$, which implies that $\mathbf u>0$ and hence $u_h>0$ in $\Omega$.
	
	In the following, we will show by contradiction that any discrete ground state satisfies $|u_h|>0$ in~$\Omega$, i.e., there cannot be a change of sign. Consider a discrete ground state $u_h = \sum_{j=1}^mu_j\Lambda_j$ with a change of sign, i.e., there exists an index pair $\{k,l\}$ such that $u_k<0$ and $u_l>0$. We denote $\bar u_h = \sum_{j=1}^m |u_j|\Lambda_j$ and recall  that it holds $|u_j|>0$ for all $j \in \{1,\dots,m\}$, which can be shown by applying the arguments from above to $\bar u_h$. Then by the irreducibility of $\mathbf S$, there exists a path $(k = i_1,\dots,i_p = l)$ with $p> 1$ of non-repeating indices connecting~$k$ and~$l$ such that  $S_{i_{q}i_{q+1}}\neq 0$ holds for all $q \in \{1,\dots,p-1\}$. There must be at least one change of sign in this patch, i.e., there exists $r \in \{1,\dots,p-1\}$ such that $u_{i_r}<0<u_{i_{r+1}}$ holds.  Note that $u_{i_r}$ or $u_{i_{r+1}}$ cannot be zero, since it holds that $|u_j|>0$ for all $j \in \{1,\dots,m\}$.  Therefore, by the M-matrix property of~$\mathbf S$ it holds that
	\begin{equation*}
		S_{i_{r}i_{r+1}}|u_{i_r}||u_{i_{r+1}}| <0<S_{i_{r}i_{r+1}}u_{i_r}u_{i_{r+1}},
	\end{equation*}  
which yields that  $\mathcal{E}_h(\bar{u}_h) < \mathcal{E}_h(u_h)$. This contradiction to the assumption that~$u_h$ is a ground state  proves that $|u_h|>0$ must hold in $\Omega$.
	
	\emph{Step 2:}
	Next, we prove the uniqueness of the discrete ground state by expressing the coordinate vector of the discrete ground state as the solution of a strictly convex minimization problem. 
To overcome the non-uniqueness caused by the sign, this minimization problem will seek the component-wise square of the ground state's coordinate vector. 
	Defining for any $\mathbf p = (p_j)_{j=1}^{k} \in \mathbb R^{k}$ with $\mathbf p\geq 0$ the norm
	\begin{equation*}
		|\mathbf{p}|_C \coloneqq \sum_{i = 1}^{\#\mathcal{T}_h}\sum_{j =1}^{d+1}\frac{|K|}{d+1} p_{(i-1)(d+1)+j},
	\end{equation*}
the desired minimization problem is then posed on the convex set
	\begin{equation*}
		C:= \{\mathbf w\in\mathbb{R}^{m}\with \mathbf w\geq 0,\;|\mathbf P \mathbf w|_{C}=1\}.
	\end{equation*}
	and seeks 
	\begin{equation}
		\label{eq:strictconvmin}
		\mathbf{v}\in \argmin_{\mathbf w\in C}\tfrac{1}{2}\sqrt{\mathbf w}^{T}\mathbf S\sqrt{\mathbf w} + \tfrac{1}{2}|\mathbf V\circ (\mathbf P\mathbf w)|_C +\tfrac{\kappa}{4}|\mathbf P \mathbf w^2|_C,
	\end{equation}
where $\sqrt{\cdot}$ denotes the component-wise square root and $\circ$ the component-wise multiplication. To show that this minimization problem is strictly convex, it suffices to verify the convexity of the first and second summands in \cref{eq:strictconvmin} and the strict convexity of the last summand (recall that $\kappa>0$). For proving the convexity of the first summand, we consider arbitrary $\mathbf v, \mathbf w \in C$ and $0\leq t\leq 1$. By the Cauchy-Schwarz inequality, it holds for all $i,j\in \{1,\dots,m\}$ that 
\begin{equation}
	\label{eq:csu}
	 t\sqrt{v_iv_j} + (1-t)\sqrt{w_iw_j} \leq \sqrt{tv_i +(1-t)w_i}\sqrt{\smash[b]{
	 	tv_j+(1-t)w_j}}.
\end{equation}
The M-matrix property of $\mathbf S$ then implies that its off-diagonal entries are non-positive, which together with \cref{eq:csu} yields that
\begin{align*}
	&\sqrt{t\mathbf v +(1-t)\mathbf w}^{T}\mathbf S\sqrt{t\mathbf v+(1-t)\mathbf w}\\[1.25ex]
	&\qquad = \sum_{j}S_{jj}(tv_j+(1-t)w_j) + \sum_{i\neq j}S_{ij}\sqrt{tv_i +(1-t)w_i}\sqrt{\smash[b]{tv_j+(1-t)w_j}}\\[-.5ex]
	&\qquad \leq t\sqrt{\mathbf v}^{T}\mathbf S\sqrt{\mathbf v} + (1-t)\sqrt{\mathbf w}^{T}\mathbf S\sqrt{\mathbf w},
\end{align*}
which proves the convexity of the first summand. The convexity and strict convexity of the second and third summand in \cref{eq:strictconvmin}, respectively, follows immediately. The unique existence of a solution $\mathbf v$ to \cref{eq:strictconvmin} then follows by classical convex optimization theory. Noting that $\pm u_h$ minimizes \cref{eq:gsdisc} if and only if $\mathbf v = \mathbf u^2$ minimizes \cref{eq:strictconvmin}, the unique existence up to sign of the discrete ground state can be concluded. 
\end{proof}

Note that in the following we will always choose the signs of the ground state~$u$ and its discrete approximation $u_h$ so that  both functions are positive in~$\Omega$. The next theorem, which is the second major result of this paper, shows that the discrete ground state eigenvalue $\lambda_h$ is the smallest eigenvalue of the nonlinear eigenvalue problem~\cref{eq:discrete_eigen} and that $\lambda_h$ is a simple eigenvalue. We emphasize that until now mainly the properties of the linearized discrete Gross-Pitaevskii eigenvalue problem have been studied (see, e.g., \cite{positivitygpe}), while the properties of the nonlinear discrete eigenvalue problem are not well understood. 
\begin{theorem}[Discrete ground state eigenvalue]\label{theo:eigenvalue:mini_simp}
	Suppose that the stiffness matrix~$\mathbf S$ is an irreducible M-matrix. Then for any eigenpair $(v_h,\mu_h)\in V_h^0\times \mathbb{R}$ of~\cref{eq:discrete_eigen}, it holds either that $\mu_h> \lambda_h$ or that $\mu_h =\lambda_h$ and $v_h = \pm u_h$. Therefore, the ground state eigenvalue $\lambda_h$ is the smallest eigenvalue of~\cref{eq:discrete_eigen} and it is simple.  
\end{theorem}
\begin{proof}
	The proof is again done in two steps: In Step 1, we consider a linearized version of the discrete eigenvalue problem~\cref{eq:discrete_eigen} and prove preliminary results, which are then used in Step 2 to prove the assertion. The following proof is inspired by~\cite[Lem.~2]{CCM10}, where a similar result is proved in the continuous~setting.
	
	\emph{Step 1:}
		To derive a linearized version of~\cref{eq:discrete_eigen} we freeze the nonlinearity in the discrete ground state~$u_h$. The resulting linearized eigenvalue problem seeks eigenpairs $(v_h,\mu_h)\in V_h^0\times \mathbb{R}$ with $\|v_h\|_{\ell}=1$ such that
\begin{equation}\label{eq:eigen:linear}
	(\nabla v_h,\nabla w_h)_{L^2} + \ell(Vv_h,w_h) + \kappa \ell(|u_h|^2v_h,w_h) = \mu_h\ell(v_h,w_h)
\end{equation}
holds for all $w_h\in V_h^0$. We note that the discrete ground state eigenpair $(u_h,\lambda_h)$ is also an eigenpair of problem \cref{eq:eigen:linear}. 
	By the min-max principle, the smallest eigenvalue of \cref{eq:eigen:linear}, denoted by $\mu_{h,1}$, can be characterized as 
\begin{equation*}\label{eq:mini:linear}
	\mu_{h,1} = \inf_{v_h \in V_h^0\with \|v_h\|_\ell = 1}(\nabla v_h,\nabla v_h)_{L^2} + \ell(Vv_h,v_h) + \kappa \ell(|u_h|^2v_h,v_h),
\end{equation*}
and the associated eigenstate, denoted by $v_{h,1}$, is the state where the minimum is attained. Using the same arguments as in the proof of \cref{theo:uni_posi:mini}, one can prove that $|v_{h,1}|>0$ holds in $\Omega$. As a consequence, we obtain that  $\ell(v_{h,1},u_h) \neq 0$, which in turn implies that $\mu_{h,1}= \lambda_h$ and that $\mu_{h,1}$ is a simple eigenvalue of~\cref{eq:eigen:linear}.

	\emph{Step 2:} Next, we return to the nonlinear eigenvalue problem \cref{eq:discrete_eigen}. We consider an arbitrary eigenvector $v_h$ of \eqref{eq:discrete_eigen}, which we write as the linear combination $v_h =\sum_{j=1}^m v_j\Lambda_j$. A non-negative version of $v_h$ can be defined as $\bar{v}_h = \sum_{j=1}^m|v_j|\Lambda_j$ and we denote $\tilde w_h \coloneqq  \bar{v}_h-u_h$. In the case that $\tilde{w}_h \leq 0$ in $\Omega$, one obtains with $\|\bar{v}_h\|_{\ell} = \|u_h\|_{\ell} = 1$ that $\bar{v}_h =u_h$. This implies that $(v_h,\mu_h)$ is also an eigenpair of the linearized eigenvalue problem~\eqref{eq:eigen:linear}. Therefore, using that $\lambda_h$ is the smallest eigenvalue of the linearized eigenvalue problem~\cref{eq:eigen:linear} and that it is simple, yields either that $\mu_h>\lambda_h$ or  that $v_h = \pm u_h$ and $\mu_h = \lambda_h$. 
	In all other cases, there exists a node $p$ of the mesh $\mathcal{T}_h$ with $\tilde{w}_h(p) >0$. We consider the function ${w_h\coloneqq v_h-u_h = \sum_{i=1}^mw_i\Lambda_i}$ which, after possibly replacing~$v_h$ by~$-v_h$, satisfies  $w_h(p)>0$. We can split the function $w_h$ as $w_h = w_h^+ + w_h^-$, where $w_h^+ = \sum_{j=1}^mw^{+}_j\Lambda_j$ with $w^+_{j} =\max(w_j, 0)\geq 0$ and $w^- = \sum_{j=1}^mw^{-}_j\Lambda_j$ with $w^-_{j} = \min (w_j,0)\leq 0$. Testing the eigenvalue problems for $v_h$ and $u_h$, cf.~\cref{eq:discrete_eigen}, with~$w_h^+$ and subtracting the resulting equations yields that
\begin{align}\label{eq:diffevp}
	\begin{split}
		&(\nabla w_h,\nabla w_h^+)_{L^2} + \ell(Vw_h,w_h^+) +  \kappa \ell(u_h^2w_h,w_h^+) -\lambda_h\ell(w_h,w_h^+)\\
		&\qquad\quad\;\; + \kappa\ell((v_h^2-u_h^2)v_h,w_h^+)
		= (\mu_h-\lambda_h)\ell(v_h,w_h^+).
	\end{split}
\end{align}
In the following, we prove that the left-hand side of \cref{eq:diffevp} is positive. Noting that, by the M-matrix property of $\mathbf S$, it holds that
\begin{equation*}
	(\nabla w_h,\nabla w_h^+)_{L^2} = (\nabla (w_h^+ + w^-),\nabla w_h^+)_{L^2}\geq (\nabla w_h^+ ,\nabla w_h^+)_{L^2}
\end{equation*}
and that
\begin{align*}
	\ell(Vw_h^-,w_h^+) = \ell(u_h^2w_h^-,w_h^+) = \ell(w_h^-,w_h^+) = 0
\end{align*}
by definition \cref{eq:lumpedbilinearform} of the lumped bilinear form, we obtain the estimate
\begin{align*}
	&(\nabla w_h,\nabla w_h^+)_{L^2} + \ell(Vw_h,w_h^+) +  \kappa \ell(u_h^2w_h,w_h^+)-\lambda_h\ell(w_h,w_h^+)\\
	&\qquad \geq (\nabla w_h^+,\nabla w_h^+)_{L^2} + \ell(Vw_h^+,w_h^+) +  \kappa \ell(u_h^2w_h^+,w_h^+)-\lambda_h\ell(w_h^+,w_h^+)\geq 0.
\end{align*}
Together with the estimate
\begin{equation*}
	\ell((v_h^2-u_h^2)v_h,w_h^+) = \ell(v_h(v_h+u_h)w_h,w_h^+) = \ell(v_h(v_h+u_h)w_h^+,w_h^+)>0
\end{equation*}
and recalling that $\kappa>0$, we obtain the positivity of the left-hand side of \cref{eq:diffevp}. Since  $\ell(v_h,w_h^+) >0$, it must holds that $\mu_h>\lambda_h$ which concludes the proof.
\end{proof}

\section{Global convergence to discrete ground state}\label{sec:globconv}

In this section, we present a fully discretized Sobolev gradient flow and prove its global convergence to the discrete ground state. To prove the global convergence, one needs to identify the limit of the fully discretized gradient flow, which can be done using the following theorem. The theorem proves that any non-negative discrete eigenstate of~\cref{eq:discrete_eigen} must necessarily coincide with the discrete ground state. The proof of this result is based on a discrete version of Picone's inequality, cf.~\cite{Brasco2014}. Note that such an inequality is also key to proving a similar result in the continuous setting; see, e.g.,~\cite[Lem.~5.4]{HenP20}.

\begin{theorem}[Non-negative discrete eigenstates]\label{thmposexcstate}
	Suppose that the stiffness matrix~$\mathbf S$ is an irreducible M-matrix and let $(v_h,\mu_h) \in V_h^0 \times R$ be an eigenpair of~\cref{eq:discrete_eigen}. Then, if $v_h\geq0$, it must hold that $v_h =  u_h$ and $\mu_h = \lambda_h$. Therefore, any non-negative discrete eigenstate must coincide with the discrete ground state.
\end{theorem}
\begin{proof}
	Also this proof is done in two steps: In Step 1, we first prove a discrete version of Picone's identity, which is then used in Step 2 to conclude the assertion. 
	
	\emph{Step 1:}
	Let us consider arbitrary vectors $\mathbf u, \mathbf v \in \mathbb R^m$ with $\mathbf u \geq 0$ and $\mathbf v >0$. One can prove by Young's inequality that for the components of $\mathbf u$ and $\mathbf v$, denoted by~$u_j$ and $v_j$, respectively, it holds for all $i,j\in \{1,\dots,m\}$ that
	\begin{equation}
		\label{eq:young}
		u_iu_j \leq \frac12 \frac{u_i^2}{v_i}v_j + \frac12 \frac{u_j^2}{v_j}v_i. 
	\end{equation}
	The matrix $\mathbf S$ is a symmetric positive definite M-matrix, which implies that its diagonal entries are positive and its off-diagonal entries are non-positive. Interpreting the square and the division of vectors component-wise, the symmetry of~$\mathbf S$ and~\cref{eq:young} yield that
	\begin{align}\label{eq:picone}
		\begin{split}
			\langle \mathbf S \mathbf v,\mathbf u^2/\mathbf v\rangle &= \sum_{i,j} S_{ij} \frac{u_i^2}{v_i}v_j = \sum_j S_{jj} u_j^2 + \sum_{j \neq i}S_{ij}\bigg(\frac12 \frac{u_i^2}{v_i}v_j + \frac12 \frac{u_j^2}{v_j}v_i \bigg)\\
			&\leq \sum_{i,j} S_{ij} u_iu_j = \langle \mathbf S \mathbf u,\mathbf u \rangle,
		\end{split}
	\end{align} 
	where $\langle\cdot ,\cdot \rangle$ denotes the Euclidean inner product of $\mathbb R^m$. 
	This inequality  can be considered a discrete version of Picone's inequity.
	
	\emph{Step 2:}
	First, note that using the arguments from the proof of \cref{theo:uni_posi:mini}, one can prove that $v_h\geq 0$ actually implies that $v_h>0$ holds in $\Omega$. Therefore,~$u_h$ and~$v_h$ are both positive discrete eigenstates of \cref{eq:discrete_eigen} with the eigenvalues~$\lambda_h$ and~$\mu_h$, respectively.  Due to their positivity, we can define the test function $w_h \coloneqq I_h(u_h^2/v_h)$, where we set $w_h$ to zero for all boundary nodes, i.e., $w_h  \in V_h^0$. Here $I_h\colon \mathcal C^0(\overline{\Omega})\to V_h$ denotes the nodal interpolation. Note that, by the normalization condition $\|u_h\|_\ell = 1$, it also holds that $\ell(v_h,w_h) = 1$. Applying the  discrete Picone inequality, cf.~\cref{eq:picone}, for the coordinate vectors  $\mathbf u$ and $\mathbf v$  of the representation of~$u_h$ and~$v_h$ in terms of the hat functions,  we obtain that
	\begin{align*}
		\mu_h &= \mu_h \ell (v_h,w_h) = (\nabla v_h,\nabla w_h)_{L^2} + \ell(Vv_h,w_h) + \kappa  \ell(|v_h|^2v_h,w_h)\\
		&\leq (\nabla u_h,\nabla u_h)_{L^2}  + \ell(Vu_h,u_h) + \kappa \ell(|u_h|^2,|v_h|^2)\\
		& \leq \lambda_h - \tfrac{\kappa}{2}\ell(|u_h|^2,|u_h|^2) + \tfrac{\kappa}{2} \ell (|v_h|^2,|v_h|^2).
	\end{align*}
	In the last step we have used~\cref{eq:lumpedbilinearform} and Young's inequality. We conclude that
	\begin{align*}
		2 \mathcal E_h (v_h) = \mu_h  -\tfrac{\kappa}{2}\ell(|v_h|^2,|v_h|^2) \leq \lambda_h -\tfrac{\kappa}{2}\ell(|u_h|^2,|u_h|^2) = 2 \mathcal E_h(u_h).
	\end{align*}
	Therefore, it must hold that $v_h$ is a ground state. Due to the uniqueness of the discrete ground state, cf.~\cref{theo:uni_posi:mini}, it follows that $v_h = u_h$ and $\mu_h = \lambda_h$.
\end{proof}

In the following, we present a method for solving the discrete constrained minimization problem \cref{eq:gsdisc} in practice. The method is an application of the energy-adaptive Sobolev gradient flow~\cite{HenP20} to the present discrete setting. There it is proved that the iteration, obtained after discretizing the gradient flow, converges globally to the continuous ground state. In this section, we will prove that a similar convergence result also holds after discretization in space, i.e., the resulting fully discretized gradient flow converges globally to the discrete ground state.

To define the fully discretized gradient flow, we introduce the discrete Green's operator of the Gross-Pitaevskii problem. This operator is henceforth denoted by $\mathcal G^h_{w_h}\colon V_h \to V_h^0$ and, for a fixed $w_h \in V_h^0$, is defined as the map of a source term~$f_h$ to the solution $u_h$, which is uniquely defined by setting for all $v_h \in V_h^0$
\begin{equation}
	(\nabla u_h,\nabla v_h)_{L^2} + \ell(Vu_h,v_h) + \ell(|w_h|^2u_h,v_h) = \ell(f_h,v_h).
\end{equation} 
Given the initial guess $u_h^0 \in V_h^0$ with $\|u_h^0\|_\ell = 1$, the iterates of the fully discretized gradient flow are for all $n = 0,1,2,\dots$ defined as  
\begin{align}
	\label{eq:discgradflow}
	\tilde u_h^{n+1} = (1-\tau^n) u_h^n + \tau^n (u_h^n,\mathcal G_{u_h^n}^h u_h^n)_{L^2}^{-1}\mathcal G_{u_h^n}^h u_h^n,\qquad u_h^n \coloneqq \frac{\tilde u_h^n}{\|\tilde u_h^n\|_{\ell}},
\end{align}
where $(\tau^n)_{n = 0}^\infty$ is a sequence of positive step sizes. 

To prove the global convergence of this iteration, we need to impose two conditions on the chosen sequence of step sizes.  First, in order to prevent stagnation of the iteration, we require that there exists a constant $c>0$ such that
\begin{equation}
	\label{eq:rest1}
	\tau^n \geq c
\end{equation} 
holds for all $n$. Second, we need to impose an upper bound on the sequence of step sizes, which can be derived following the lines of \cite[Lem.~4.7]{HenP20}. It reads
\begin{equation}
	\label{eq:rest2}
	\tau^n \leq 2\min\{(1+\kappa C_1C_2^4)^{-1},\mathcal E_h(u_h^0)^{-1/2}\}
\end{equation}
for all $n$, where $C_1, C_2>0$ are constants that can be bounded by explicitly computable expressions. More precisely, the constant $C_1$ is the norm equivalence constant satisfying $\|v_h^2\|_{\ell}^2\leq C_1\|v_h\|_{L^4}^4$ for all $v_h \in V_h^0$. Using a transformation to the reference simplex, one derives the explicit bound 
\begin{equation*}
	C_1 \leq \frac{1}{\gamma (d+1)d!},
\end{equation*}
where $\gamma>0$ denotes the smallest eigenvalue of the element mass matrix corresponding to the reference simplex for the quadratic FEM (using the Lagrange basis). The constant $C_2$ is the continuity constant of the embedding $H^1_0(\Omega)\hookrightarrow L^4(\Omega)$. Using Hölder's inequality and the Gagliardo-Nirenberg-Sobolev inequality (see, e.g., \cite[Sec.~5.6.1]{Eva10}), we obtain the explicit  bound
\begin{equation*}
	C_2 \leq \left\{
	\begin{array}{lr}
		|\Omega|^{3/4} & d = 1,\\
		2|\Omega|^{1/4} & d = 2,\\
		4|\Omega|^{1/12} & d = 3,
	\end{array}\right.
\end{equation*}
where $|\Omega|$ denotes the volume of $\Omega$. 
Inserting the above  bounds into \cref{eq:rest2} gives an upper bound for the step sizes which is explicit in $\kappa,d,\gamma,\mathcal E(u_h^0),$ and $|\Omega|$. 

The following corollary proves a global convergence result for the fully discretized gradient flow~\cref{eq:discgradflow}. 
For the global convergence proof, \cref{thmposexcstate} is of great importance, since it allows to identify the (non-negative) limit function of iteration~\cref{eq:discgradflow} as the discrete ground state. To the best of our knowledge, this is the first global convergence result to the discrete ground state in the fully discrete setting. Our arguments are not restricted to the fully discretized gradient flow~\cref{eq:discgradflow}. In fact, the arguments apply to any iteration used for the numerical solution of~\cref{eq:gsdisc} that converges globally to a stationary state and preserves the non-negativity of the iterates, e.g., the $J$-method from \cite{AltHP21} with a suitable shift or the $H^1$-gradient flow from \cite{positivitygpe}.
Note that it is  possible to quantify the rate of local linear convergence of the above methods, along with a corresponding convergence radius. This has be done, e.g., for the gradient flow of~\cite{HenP20} with an energy-adaptive metric in \cite{Zha21ppt,AltmannPeterseimStykel2022,MR4595654}, for the $H^1$-gradient flow in \cite{DaK10,positivitygpe}, and for the damped $J$-method in~\cite{AltHP21}.

\begin{corollary}[Global convergence to discrete ground state]\label{cor:globconv}
	Suppose that the stiffness matrix~$\mathbf S$ is an irreducible M-matrix and let the step size condition $\tau^n \leq 1$ for all $n$ as well as conditions  \cref{eq:rest1,eq:rest2} be satisfied. Then given a non-negative  initial guess $u_h^0 \geq 0$, which is normalized with respect to the lumped $L^2$-norm, the sequence $(u_h^n)_{n = 0}^\infty$ defined by~\cref{eq:discgradflow} satisfies for all $n$ that $u_h^n \geq 0$. Furthermore, the sequence of iterates converges to the ground state $u_h$ defined in \cref{eq:gsdisc}.	
\end{corollary}

\begin{proof}
	Similar to \cite[Cor.~4.11]{HenP20}, under the step size conditions \cref{eq:rest1,eq:rest2}, it can be shown that the energy is strictly reduced and that there exists a limit energy $E_h^* := \lim_{n\rightarrow\infty}\mathcal{E}_h(u_h^n)$. Similar to the existence proof of minimizers for the finite-dimensional minimization problem~\cref{eq:gsdisc}, up to subsequences, we obtain that $u_h^n\rightarrow v_h$ where $v_h\in V_h^0$ with $\|v_h\|_{\ell}=1$ and $\mathcal{E}_h(v_h) = E_h^{*}$ is a discrete eigenstate. Note that specifying the norm in which we have convergence is not important, since all norms are equivalent in the finite-dimensional setting. Furthermore, it follows from the definition of the iteration, cf.~\cref{eq:discgradflow}, that for step sizes $\tau_n \leq 1$ the iteration preserves the non-negativity of the initial iterate. As a consequence, the limit eigenstate $v_h$ is also non-negative. By \cref{thmposexcstate} it must hold that $v_h = u_h$, where~$u_h$ denotes the discrete ground state defined in~\cref{eq:gsdisc}. Note that, in the following, we will write $u_h$ instead of $v_h$ for the limit of iteration~\cref{eq:discgradflow}. 
	
	It remains to show that the whole sequence $(u_h^n)_{n = 0}^\infty$ converges to $u_h$. For this we will use that~$(u_h,\lambda_h)$ is the eigenpair corresponding to the smallest eigenvalue of the linearized eigenvalue problem~\cref{eq:eigen:linear}, which was shown in Step~1 of the proof of \cref{theo:eigenvalue:mini_simp}, and that $u_h$ is positive in $\Omega$ according to \cref{theo:uni_posi:mini}. Following the arguments of \cite[Thm.~5.1]{HenP20}, one then obtains that
	\begin{equation*}
		\ell(|u_h^n|^2-|u_h|^2,|u_h^n|^2-|u_h|^2) \leq \tfrac{4}{\kappa}(\mathcal{E}_h(u_h^n)-E_h^*)\rightarrow 0.
	\end{equation*} 
From the latter convergence result one can conclude that $|u_h^n|^2 \to |u_h|^2$, which, due to the non-negativity of $u_h^n$ and $u_h$, implies that $u_h^n \to u_h$, i.e., the convergence of the whole sequence to~$u_h$. This completes the proof.	
\end{proof}

\section{A priori error analysis}\label{sec:erroranalysis}

This section performs an a priori error analysis of the proposed mass-lumped finite element discretization for the Gross-Pitaevskii problem.  Recall that we choose the signs of~$u$ and~$u_h$ such that $u,u_h>0$ holds in $\Omega$.
To simplify the notation below, we abbreviate the ground state energy and its discrete counterpart by 
\begin{equation*}
	E \coloneqq \mathcal E(u),\qquad E_h \coloneqq \mathcal E_h(u_h)
\end{equation*}
and introduce a notation that hides constants independent of $h$ in estimates.
	\begin{remark}[Tilde notation]
If it holds that $a \leq C b$, where $C>0$ is a constant that may depend on the domain, the mesh regularity, the coefficients $V$ and $\kappa$, and the ground state $u$, but is independent of the mesh size $h$, we may write $a \lesssim b$  to hide the constant. Analogously, we may write $b\gtrsim a$ for $a \geq C b$.
\end{remark}

The following theorem proves optimal orders of convergence for the ground state, energy, and eigenvalue approximations of the proposed method. 

\begin{theorem}[A priori error analysis]\label{theo:plain_convergence}
	Assume that $V$ is $\mathcal T_h$-piecewise $H^2$-regular with an uniformly bounded piecewise $H^2$-norm. Then, the ground state approximations $u_h$ defined in~\cref{eq:gsdisc} converge to the  ground state $u$ defined in~\cref{eq:gs} with 
	\begin{equation}\label{eq:errl2h1}
		\|u-u_h\|_{H^1}\lesssim h,\qquad \|u-u_h\|_{L^2} \lesssim h^2.
	\end{equation} 
Furthermore, the energies and eigenvalue approximations $E_h$ and $\lambda_h$ converge to their continuous counterparts $E$ and $\lambda$, respectively, with 
\begin{equation}
	\label{eq:convenergylambda}
	|E-E_h|\lesssim h^2,\qquad |\lambda-\lambda_h|\lesssim h^2.
\end{equation}
\end{theorem}

\begin{proof}
	This proof is done in two steps. In Step 1 we prove the boundedness of the discrete energies $E_h$ and their second-order convergence to $E$. This result is then used to establish the first-order $H^1$-convergence of the ground state approximations~$u_h$ to the ground state~$u$. In Step 2 we employ  a duality argument to prove the second-order $L^2$-convergence of $u_h$ to $u$, which also implies the second-order convergence of the eigenvalue approximations~$\lambda_h$ to~$\lambda$. 
	
	\emph{Step 1:}
	In this proof we utilize~\cite[Thm.~3]{CCM10}, which proves the convergence of the standard linear FEM to the 
	ground state. The ground state approximations of this method are henceforth denoted by $\hat u_h \in V_h^0$. Using \cref{lem:estlump} we get that
	\begin{align}
		\label{eq:wrongnormest}
		|\|\hat u_h\|_\ell-1|(\|\hat u_h\|_\ell+1) = |\|\hat u_h\|_\ell^2-1| \lesssim h^2\|\nabla \hat u_h\|_{L^2}^2 \lesssim h^2,
	\end{align}
which yields the estimate $|\|\hat u_h\|_\ell-1| \lesssim h^2$. This estimate allows us to estimate the error between $\hat u_h$ and its rescaled version $\tilde{u}_h \coloneqq \hat{u}_h/\|\hat{u}_h\|_{\ell}$, which is normalized with respect to $\|\cdot\|_\ell$. To estimate the difference between the energies $\mathcal E_h(\tilde u_h)$ and~$E$, we employ the triangle inequality which yields that 
\begin{align*}
	\begin{split}
			|\mathcal E_h(\tilde u_h)-E| &\leq |\mathcal E_h(\tilde u_h) - \mathcal E_h(\hat u_h)| + |\mathcal E_h(\hat u_h) - \mathcal E(\hat u_h)| + |\mathcal E(\hat u_h) - E|\\
		& \eqqcolon \Xi_1+\Xi_2+\Xi_3.
	\end{split}
\end{align*}
Below the terms $\Xi_1$, $\Xi_2$, and $\Xi_3$ are estimated individually. Using \cref{eq:wrongnormest}, we obtain for the first term that $\Xi_1 \lesssim h^2$. To estimate $\Xi_2$, we use \cref{lem:estlump}, the uniformly bounded piecewise $H^2$-norm of $V$, and the uniform $L^\infty$-bound for $\hat u_h$ from \cref{lem:linftybounds} to get that
	\begin{align*}
		\Xi_2 &\leq |\ell(V\hat u_h,\hat u_h)-(V\hat u_h,\hat u_h)_{L^2}| + |\ell(|\hat u_h|^2\hat u_h,\hat u_h)-(|\hat u_h|^2\hat u_h,\hat u_h)_{L^2}|\\
		&\lesssim h^2 \big(\|\nabla \hat u_h\|_{L^2}\|\hat u_h\|_{L^\infty} + \|\hat u_h\|_{H^1}^2+ \|\hat u_h\|_{L^\infty}^2\|\nabla \hat u_h\|_{L^2}^2\big)\lesssim h^2.
	\end{align*}
The estimate $\Xi_3 \lesssim h^2$ for the third term can be found in \cite[Thm.~3]{CCM10}. Combining the above estimates yields that
\begin{align}
	\label{eq:diff}
	|\mathcal E_h(\tilde u_h)-E|\lesssim h^2.
\end{align}
Using that $E_h\leq \mathcal{E}_h(\tilde{u}_h)$ and \cref{eq:diff}, we obtain that the discrete energies $E_h$ are uniformly bounded, and hence the same applies to $\|u_h\|_{H^1}$, $\|u_h^2\|_\ell$, $\|u_h\|_{L^q}$ for all $1 \leq q \leq 6$, and~$\lambda_h$. Using \cref{lem:linftybounds}, we additionally get the uniform boundedness of $\|u_h\|_{L^{\infty}}$. Next, we define $\bar{u}_h = u_h/\|u_h\|_{L^2}$, which is a $L^2$-normalized version of~$u_h$. Similar to \cref{eq:wrongnormest}, one can prove that $|\|u_h\|_{L^2}-1| \lesssim h^2$, which yields that 
\begin{equation}
	\label{eq:estscalendminnotscaled}
	\|u_h-\bar{u}_h\|_{H^1}\lesssim h^2.
\end{equation}
To estimate the difference between the energies $E_h$ and $\mathcal E(\bar u_h)$, we employ the triangle inequality to get that 
\begin{align*}
	|E_h-\mathcal{E}(\bar{u}_h)|\leq |E_h - \mathcal E(u_h)| + |\mathcal E(u_h) - \mathcal E(\bar u_h)|,
\end{align*}
where the first term can be estimated similarly as $\Xi_2$ and the second term can be estimated using that  $|\|u_h\|_{L^2}-1| \lesssim h^2$. We obtain that
\begin{align}\label{eq:estdiffeh-ehbar}
	|E_h-\mathcal{E}(\bar{u}_h)|\lesssim h^2.
\end{align}

To derive the second-order error estimate for the discrete energies, we bound the difference $E_h -E$ from above and below. We derive the bounds
\begin{align}
	\label{eq:energylowandupbound}
	\begin{split}
			E_h-E &= E_h-\mathcal{E}_h(\tilde{u}_h) + \mathcal{E}_h(\tilde{u}_h)-E\leq \mathcal{E}_h(\tilde{u}_h)-E \lesssim h^2,\\
		E_h-E &= E_h-\mathcal{E}(\bar{u}_h) + \mathcal{E}(\bar{u}_h)-E\geq  E_h-\mathcal{E}(\bar{u}_h)\gtrsim -h^2,
	\end{split}
\end{align}
where we used~\cref{eq:diff} and $E_h \leq \mathcal{E}_h(\tilde{u}_h)$ for the upper bound as well as \cref{eq:estdiffeh-ehbar} and ${E \leq \mathcal E(\bar u_h)}$ for the lower bound. This proves the second-order convergence of the energy approximations, which is the first estimate in \cref{eq:convenergylambda}. 

To prove the first-order $H^1$-convergence of the ground state approximations, we use the triangle inequality and \cref{eq:estdiffeh-ehbar,eq:energylowandupbound} to get that $|E-\mathcal{E}(\bar{u}_h)| \lesssim h^2$. Recalling that $u_h$ and $u$ have the same signs, \cite[Thm.~1]{CCM10}  proves the first-order estimate $\|\bar{u}_h-u\|_{H^1}\lesssim h$ for the rescaled approximations $\bar u_h$. Using the triangle inequality and \cref{eq:estscalendminnotscaled}, we can conclude the first-order $H^1$-convergence of the ground state approximations, which is the first estimate of \cref{eq:errl2h1}. 

\emph{Step 2:}
Next we use a duality argument to prove the second-order $L^2$-convergence of the ground state approximations. This argument is based on the auxiliary problem of 
\cite[Eq.~(70)]{CCM10} which, for a given $w\in L^2(\Omega)$,  seeks $z\in H^1_0(\Omega)$ such~that 
\begin{equation}
	\label{eq:adjprob}
	-\Delta z + (V+3\kappa u^2 -\lambda) z
	=
	2 \kappa(u^3,z)_{L^2} u + w - (w,u)_{L^2} u
\end{equation}
holds in the weak sense. This problem is solved by the unique solution $z \in u^\perp \coloneqq \{v \in H^1_0(\Omega)\with \tspof{u}{v} = 0\}\subset H^1_0(\Omega)$ satisfying
\begin{equation}
	\label{eq:auxprob}
	J_{u,\lambda}(z,v) \coloneqq \tspof{\nabla z}{\nabla v} + \tspof{(V+3\kappa u^2-\lambda)z}{v} = \tspof{w}{v}	
\end{equation}
for all $v \in u^\perp$. The well-posedness of this problem is a consequence of the Lax-Milgram theorem using the coercivity and continuity of the bilinear form $J_{u,\lambda}$, cf.~\cite[Lem.~1]{CCM10}, and the fact that $u^\perp$ is a complete subspace of~$H^1_0(\Omega)$. Classical elliptic regularity theory on convex domains then implies that $z \in H^2(\Omega)$ with the estimate $\|z\|_{H^2} \lesssim \|w\|_{L^2}$. Note that due to the continuous embedding $H^2(\Omega)\hookrightarrow \mathcal C^0(\overline{\Omega})$ we also get the estimate $\|z\|_{L^\infty}\lesssim\|w\|_{L^2}$. 

Proceeding similarly as in the proof of \cite[Thm.~1]{CCM10}, we define the function $u_h^{*} = \bar{u}_h+ \frac{1}{2}\|\bar{u}_h-u\|_{L^2}^2u \in H^1_0(\Omega)$ and note that it holds $(u_h^{*},u)_{L^2} = 1$ since $\|\bar{u}_h\|_{L^2}=1$. Setting $w = \bar u_h - u$, using definition~\cref{eq:auxprob}, and performing a number of algebraic manipulations, we obtain that
\begin{align}
	\|w\|_{L^2}^2 
	& = (w,u_h^{*}-u)_{L^2} + \tfrac{1}{4}\|w\|_{L^2}^4 = J_{u,\lambda}(z, u_h^{*}-u) + \tfrac{1}{4}\|w\|_{L^2}^4\notag\\
	& = J_{u,\lambda}(w,z) + \tfrac{1}{2}\|w\|_{L^2}^2J_{u,\lambda}(u,z) + \tfrac{1}{4}\|w\|_{L^2}^4\notag \\
	&  = J_{u,\lambda}(w,I_hz) + J_{u,\lambda}(w,z-I_hz)+\kappa\|w\|_{L^2}^2(u^3,z)_{L^2} + \tfrac{1}{4}\|w\|_{L^2}^4,\label{eq:tobeest}
\end{align}
where $I_h\colon \mathcal C^0(\overline{\Omega}) \to V_h$ denotes the nodal interpolation. Note that, in the last step, we have used that $
	J_{u,\lambda}(u,v) = 2\kappa(u^3,v)_{L^2}$ holds for all $v \in H^1_0(\Omega)$. To estimate the terms on the right-hand side of~\cref{eq:tobeest}, we will use the following estimates 
\begin{equation}
	\label{eq:fuest}
	\|z-I_hz\|_{H^1}\lesssim h\|w\|_{L^2},\qquad \|I_h z\|_{H^1} \lesssim \|w\|_{L^2},\qquad  \|I_hz\|_{L^{\infty}}\lesssim\|w\|_{L^2},
\end{equation}
which can be proved using the properties of the nodal interpolation and the $H^2$-regularity of $z$. For the second and third terms, we obtain using these estimates~that
\begin{equation*}
	J_{u,\lambda}(w,z-I_hz)\lesssim \|w\|_{H^1}\|z-I_hz\|_{H^1}\lesssim h\|w\|_{H^1}\|w\|_{L^2}
\end{equation*}
and 
\begin{equation*}
	\kappa\|w\|_{L^2}^2(u^3,z)_{L^2}\leq \kappa\|w\|_{L^2}^2\|u^3\|_{L^2}\|z\|_{L^2}\lesssim \|w\|_{L^2}^3,
\end{equation*}
respectively.
To treat the first term on the right-hand side of~\cref{eq:tobeest}, we use that
\begin{equation*}
	(\nabla u_h,\nabla I_hz)_{L^2} =   \lambda_h\ell(u_h,I_hz)-\ell(Vu_h,I_hz) - \kappa \ell(u_h^3,I_hz),
\end{equation*}
which allows us to rearrange the term as
\begin{align*}
	J_{u,\lambda}(w,I_hz) & = J_{u,\lambda}(\bar u_h - u_h,I_hz) + J_{u,\lambda}(u_h-u,I_hz)\\
	&= J_{u,\lambda}(\bar u_h - u_h,I_hz) + \big((Vu_h,I_hz)_{L^2} - \ell(Vu_h,I_hz)\big) \\
	& \qquad + \kappa\big((u_h^3,I_hz)_{L^2}- \ell(u_h^3,I_hz)\big) + \big(\lambda_h\ell(u_h,I_hz) - \lambda(u_h,I_hz)_{L^2}\big) \\
	& \qquad  + \kappa\big(3(u^2u_h,I_hz)_{L^2}-2(u^3,I_hz)_{L^2} - (u_h^3,I_hz)_{L^2}\big)\\
	& \eqqcolon \Psi_1+\Psi_2+\Psi_3+\Psi_4+\Psi_5.
\end{align*}
Below we estimate the terms $\Psi_i$ for $i = 1,\dots,5$ individually. Using  \cref{eq:estscalendminnotscaled,eq:fuest}, we obtain for the term $\Psi_1 $ that
\begin{equation*}
	\Psi_1\lesssim\|u_h-\bar{u}_h\|_{H^1}\|I_hz\|_{H^1}\lesssim h^2\|w\|_{L^2}.
\end{equation*}
The terms $\Psi_2$ and $\Psi_3$ can be estimated with  \cref{lem:estlump} which yields that
\begin{equation*}
	\Psi_2 \lesssim h^2\|w\|_{L^2},\qquad \Psi_3 \lesssim h^2\|w\|_{L^2}.
\end{equation*}
Before considering the term $\Psi_4$, we  derive an estimate for the error of the eigenvalue approximations. Noting that $\lambda = 2E + \tfrac{\kappa}{4}\|u\|_{L^4}^4$ and $\lambda_h = 2E_h + \tfrac{\kappa}{4}\|u_h^2\|_\ell^2$, and using the first estimate of~\cref{eq:convenergylambda}, we obtain that
\begin{equation}\label{eq:err:eigenvalue}
	|\lambda-\lambda_h| \lesssim  |E-E_h| + \big|\|u\|_{L^4}^4-\|u_h\|_{L^4}^4\big| + \big|\|u_h\|_{L^4}^4-\|u_h^2\|_{\ell}^2\big|\lesssim h^2 +\|u-u_h\|_{L^2}.
\end{equation}
This result  allows us to estimate $\Psi_4$ as
\begin{align*}
	\Psi_4 &\leq \lambda_h|\ell(u_h,I_hz) -(u_h,I_hz)_{L^2}| + |(\lambda_h-\lambda)(u_h,I_hz)_{L^2}| \\
	&\lesssim h^2|I_hz|_{H^1}+|(\lambda_h-\lambda)(u_h-u,I_hz)_{L^2}|+ |(\lambda_h-\lambda)(u,z-I_hz)_{L^2}|\\
	&\lesssim h^2\|w\|_{L^2} + \big(h^2 + \|u_h-u\|_{L^2}\big)\big(\|u_h-u\|_{L^2}\|w\|_{L^2} + h\|w\|_{L^2}\big),
\end{align*}
where we used that $(u,z)_{L^2}=0$ since $z \in u^\perp$. Finally, for the term $\Psi_5$, we get that
\begin{align*}
	\Psi_5 &= |\kappa((u_h-u)^2(2u+u_h),I_hz)_{L^2}|\lesssim \|u_h-u\|_{L^2}^2\|I_hz\|_{L^{\infty}}\lesssim \|u_h-u\|_{L^2}^2\|w\|_{L^2}.
\end{align*}
Combining the above estimates then yields the following second-order estimate for the rescaled ground state approximations:
\begin{equation*}
	\|\bar{u}_h-u\|_{L^2} = \|w\|_{L^2} \lesssim h^2.
\end{equation*}
The desired $L^2$-convergence result for the ground state approximations can  be concluded using the triangle inequality and~\cref{eq:estscalendminnotscaled}. This proves the second estimate in \cref{eq:errl2h1}. The second-order convergence of the eigenvalue approximations follows directly from \cref{eq:err:eigenvalue}, which proves the second estimate in \cref{eq:convenergylambda}.
\end{proof}

\section{Numerical experiments}\label{sec:numexp}
In this section, we present numerical experiments that support the theoretical predictions of this paper. To solve the discrete minimization problem \cref{eq:gsdisc}, we use the fully discretized Sobolev gradient flow defined in~\cref{eq:discgradflow}. Note that, especially for large values of~$\kappa$, the step size bound~\cref{eq:rest2} which we needed to prove global convergence, is very restrictive (explicit values for each numerical experiment can be found in the respective subsections). Therefore, for the sake of computational efficiency, we use the adaptive choice of step sizes as outlined in \cite[Rem.~4.3]{HenP20}, where we restrict the one-dimensional minimization problem to step sizes in $[0,1]$ to ensure non-negativity of the iterates; see also~\cite{AltmannPeterseimStykel2022}. The initial iterate is constructed by interpolating a constant function to the finite element space $V_h^0$ with zero boundary conditions and normalizing the resulting function with respect to the lumped $L^2$-norm. The iteration is terminated if the relative $L^2$-residual of the current iterate falls below~$10^{-12}$.
For implementation details, see the code available at \url{https://github.com/moimmahauck/GPE_P1_lumped}, which is derived from a basic implementation of the  FEM used in~\cite{Mal20}.

 \subsubsection*{Harmonic potential with strong interaction}
\begin{figure}
	\includegraphics[height=.4\linewidth]{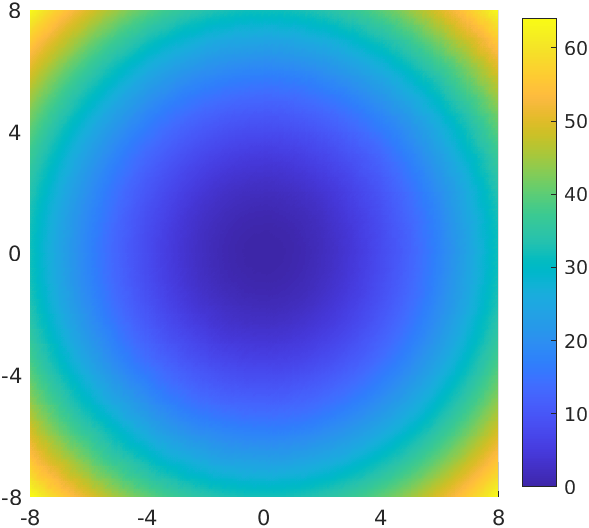}\hfil
	\includegraphics[height=.4\linewidth]{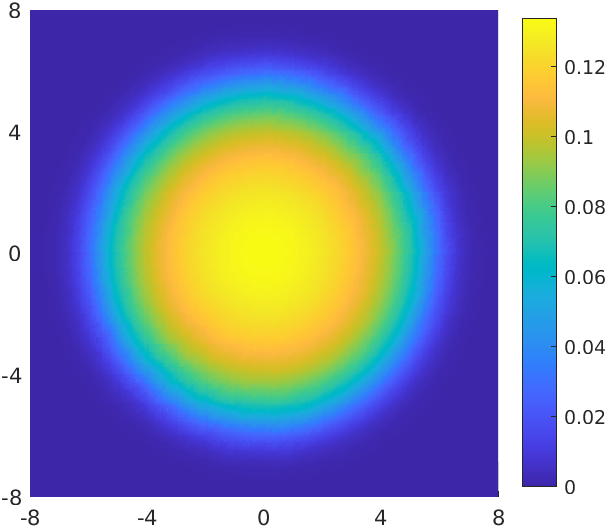}
	\caption{Illustration of the harmonic potential on the left and and a discrete ground state approximation on the right.}
	\label{fig:gsharmonic}
\end{figure}

 The problem considered in the first numerical experiment is posed on the domain $\Omega = (-8,8)^2$. We consider the harmonic potential $V(x) = \tfrac12 |x|^2$ and the particle interaction parameter $\kappa = 1000$.   The corresponding ground state is point symmetric with respect to the origin and decays exponentially. For a depiction of the harmonic potential and an approximation to the ground state, we refer to \cref{fig:gsharmonic}. Note that for this parameter setting with a large~$\kappa$, the step size bound of~\cref{eq:rest2} takes a very small value of about~$5\times 10^{-7}$. This value practically means a stagnation of the gradient flow algorithm and is therefore not feasible in practice, which explains why we use the adaptive algorithm described above.
To verify the optimal orders of convergence of the proposed method, we consider a hierarchy of Friedrichs-Keller triangulations generated by successive uniform red refinement of an initial triangulation consisting of two triangles.  Note that since no analytical solution is available, all errors are computed with respect to a reference solution. This reference solution is computed using the standard linear FEM on the mesh obtained by twice uniform red refinement of the finest mesh in the considered hierarchy. In an abuse of notation, we denote the reference ground state and the reference energy and eigenvalue by $u$, $E$, and $\lambda$, respectively. We compare the approximations of the proposed method with those of the standard linear FEM on the same mesh. For the standard linear FEM, the potential is integrated exactly using a quadrature rule of sufficiently high order. The ground state, energy, and eigenvalue approximations of the standard linear~FEM are denoted by $\hat u_h$, $\hat E_h$, and~$\hat \lambda_h$, respectively.  Note that for both spatial discretizations all iterates remain non-negative in~$\Omega$. For the proposed mass-lumped FEM this could be proved in \cref{cor:globconv}. One observes that the fully discretized gradient flow~\cref{eq:discgradflow} effectively minimizes the discrete energy from iterate to iterate.  It takes about fifty iterations for this problem to reach the specified tolerance.

\begin{figure}
 	\includegraphics[height=.45\linewidth]{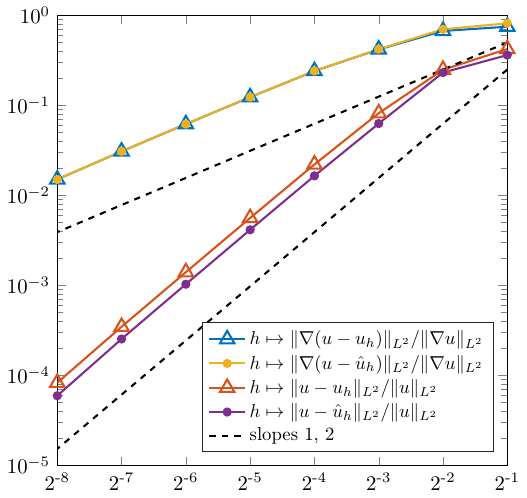}\hfil
 	\includegraphics[height=.45\linewidth]{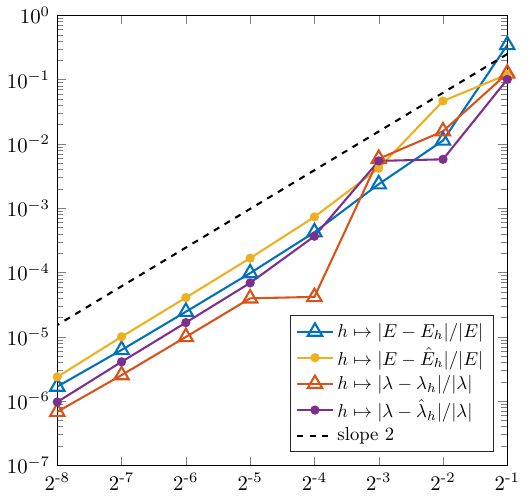}
 	\caption{Error plots for the proposed method and the  standard linear FEM for the harmonic potential. The relative  $L^2$-approximation errors of the ground state and its gradient are shown on the left. On the right, the relative energy and eigenvalue approximation errors are shown.}
 	\label{fig:convharmonic}
 \end{figure}

\begin{figure}
 	\includegraphics[height=.45\linewidth]{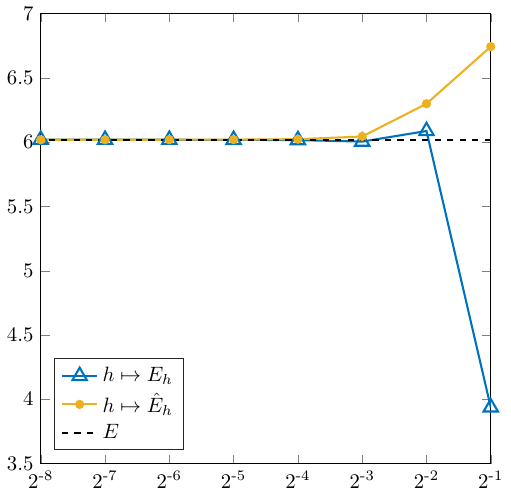}\hfil
 	\includegraphics[height=.45\linewidth]{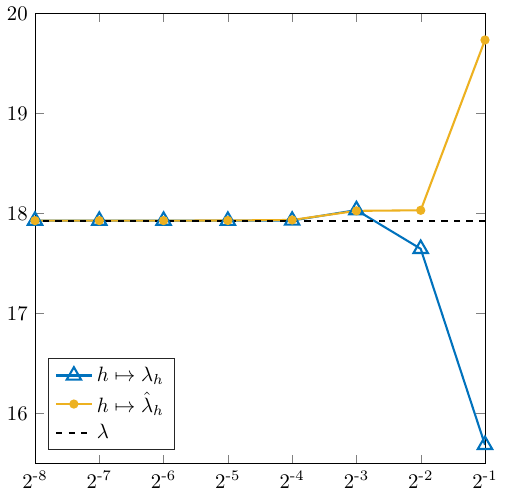}
 	\caption{Energy and eigenvalue approximations shown on the left and right, respectively, computed using the  proposed  method and the standard linear FEM.}
 	\label{fig:ergandlam}
 \end{figure}

\cref{fig:convharmonic} compares the convergence behavior of the proposed mass-lumped FEM to that of the standard linear FEM. One observes optimal convergence orders for the proposed method, which is consistent with the theoretical prediction in \cref{theo:plain_convergence}. The same convergence behavior can be observed for the standard linear FEM. Interestingly, the eigenvalue and energy approximations of the proposed lumped discretization seem to be slightly better than those of its non-lumped counterpart.
 Next, in \cref{fig:ergandlam} we examine the energy and eigenvalue approximations of the proposed method and compare them to those of the standard linear FEM. One observes that the standard linear FEM approximates the ground state energy from above due to its conformity. This is generally not true for the proposed lumped method, which can also be observed for lumped FEMs in the context of linear elliptic eigenvalue problems; see, e.g., \cite{Andreev1992}.
  
\subsubsection*{Disorder potential with exponential localization}

\begin{figure}
	\includegraphics[height=.4\linewidth]{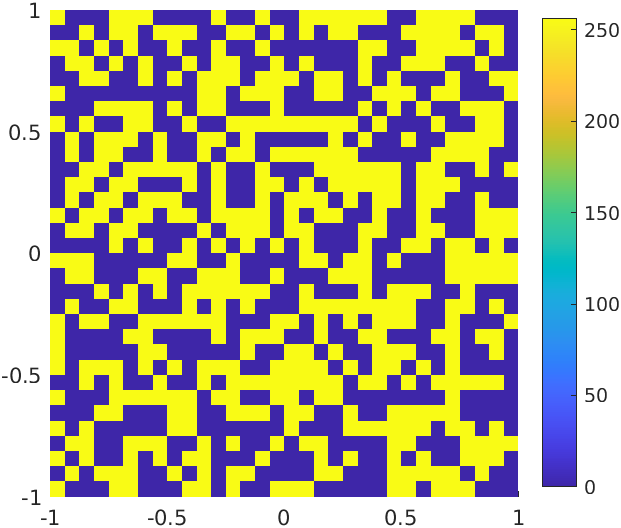}\hfil
	\includegraphics[height=.4\linewidth]{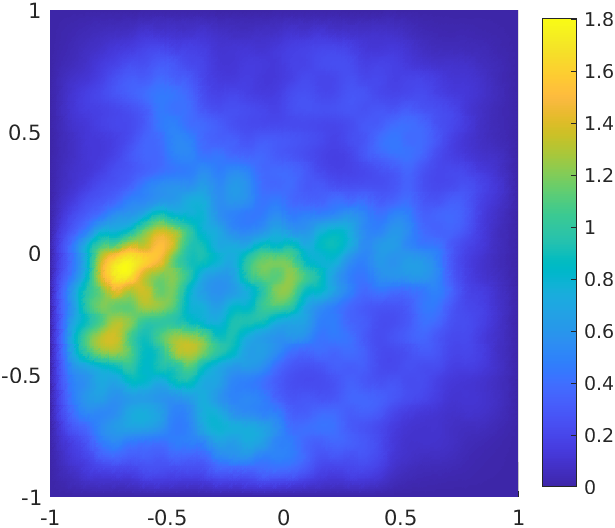}
	\caption{Illustration of the disorder potential on the left and and a discrete ground state approximation on the right.}
	\label{fig:disorder}
\end{figure}

For the second numerical experiment we consider a disorder potential on the domain $\Omega = (-1,1)^2$. This potential is constructed using a Cartesian grid of the domain with $2^5$ elements in each dimension. The potential is then chosen as the piecewise constant function on this Cartesian grid, whose element values are given by realizations of independent coin toss random variables taking the values $0$ and $(2\epsilon)^{-2}$; see \cref{fig:disorder}~(left) for an illustration. The particle interaction parameter $\kappa$ is chosen to be one. For such coefficients there occurs an effect called Anderson localization (see, e.g., \cite{AltPV18,AltHP20,AltHP22} for numerical and theoretical studies), which enforces an exponential localization of the ground state; see \cref{fig:disorder} (right). For the parameter setting in this numerical experiment, the step size bound of~\cref{eq:rest2} takes a value of about $4 \times 10^{-3}$. Due to the smaller $\kappa$, this bound is less restrictive than that of the previous numerical experiment. Nevertheless, for better computational efficiency, we use the adaptive choice of step sizes as outlines above.
Similar to the previous numerical experiment, we consider a hierarchy of Friedrichs-Keller triangulations. To compute the reference solution, we again use the standard linear FEM on the mesh obtained by twice uniform red refinement of the finest mesh in the hierarchy. 
We emphasize that this example is numerically quite challenging, as can be seen from the comparatively large number of iterations required. While the fully discretized Sobolev gradient flow method of~\cref{eq:discgradflow} required about fifty iterations for the previous numerical example, it takes several hundred iterations for this numerical example to converge to the specified tolerance. This discrepancy is related to the fact that the spectral gap to the second eigenvalue that determines the local linear rate of convergence, cf.~\cite{MR4595654}, scales with the small parameter~$\epsilon$.

\begin{figure}
	\includegraphics[height=.45\linewidth]{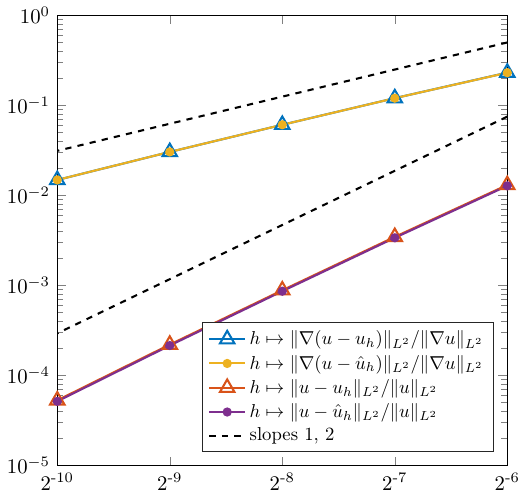}\hfil
	\includegraphics[height=.45\linewidth]{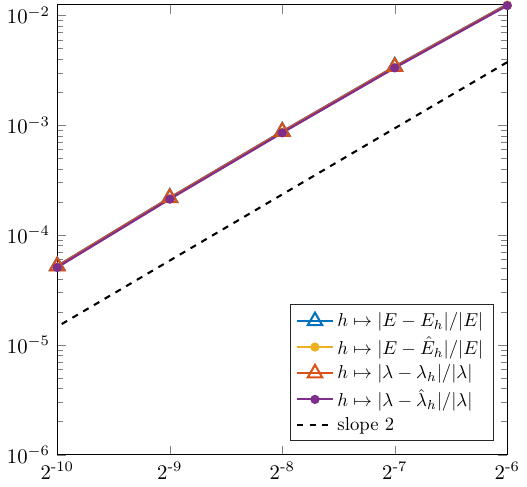}
	\caption{Error plots of the proposed method and the standard linear FEM for the disorder potential. The relative $L^2$-approximation errors of the ground state and its gradient are shown on the left. On the right, the relative energy and eigenvalue approximation errors are shown.}
	\label{fig:disorderconv}
\end{figure}

Also for this numerical example, it can be observed that the discrete ground states of the proposed method and the standard linear FEM are positive in $\Omega$. For the proposed method this was shown in~\cref{cor:globconv}. Generally, it seems difficult to construct numerical examples where the positivity is violated for the standard linear FEM. Furthermore, in \cref{fig:disorderconv}, we also observe the optimal order of convergence of both methods as the mesh size is decreased. This again supports the theoretical predictions of \cref{theo:plain_convergence}. Note that the error curves are almost on top of each other, which makes it difficult to distinguish between them.

\section{Conclusion}
In this paper, we have proposed a mass-lumped FEM for the approximation of the Gross-Pitaevskii ground state. This method is able to preserve many properties of the continuous ground state, such as positivity and uniqueness up to sign, or a Picone-type inequality. The latter paves the way for proving the global convergence of fully discretized gradient flow methods to the discrete ground state. We also prove that the proposed method has the same order of convergence as the standard linear FEM. The proposed method enjoys certain computational advantages over the standard linear FEM, e.g., a computationally cheaper assembly of the (diagonal) nonlinear term in each iteration of the fully discretized gradient flow method.

\appendix

\section{Collection of frequently used bounds}

The following lemma provides estimates for the lumping error, which are an important ingredient in the convergence proof of \cref{theo:plain_convergence}.

\begin{lemma}[Lumping error]
	\label{lem:estlump}
	Given the potential $V$, which is assumed to be $\mathcal T_h$-piecewise \mbox{$H^2$-regular} with uniformly bounded piecewise $H^2$-norm,  it holds for all $v_h,w_h\in V_h$ that 
\begin{align}
\label{eq:lumperr1}
	\begin{split}
		|\ell(Vv_h,w_h)-(Vv_h,w_h)_{L^2}| &\lesssim h^2 \big(\|\nabla v_h\|_{L^p}\|\nabla w_h\|_{L^q} + \|\nabla v_h\|_{L^r}\| w_h\|_{L^s} \\
		&\qquad\quad\, + \|v_h\|_{L^t}\|\nabla w_h\|_{L^u}+ \|v_h\|_{L^v}\|w_h\|_{L^w}\big),
	\end{split}
\end{align}
where $1\leq p,\,q,\,r,\,s,\,t,\ u,\ v,\ w \leq \infty$ are arbitrary numbers satisfying that
\begin{equation*}
	\tfrac1p+\tfrac1q = 1,\qquad  \tfrac1r+\tfrac1s=\tfrac1t+\tfrac1u= \tfrac56,\qquad  \tfrac1v+\tfrac1w=\tfrac12.
\end{equation*}
Furthermore, it holds for all $v_h,w_h\in V_h$ that
\begin{align}
	\label{eq:lumperr2}
	\begin{split}
			|\ell(|w_h|^2w_h,v_h)-(|w_h|^2w_h,v_h)_{L^2}| &\lesssim h^2 \big(\|\nabla w_h\|_{L^s}\|w_h\|_{L^t}^2\|\nabla v_h\|_{L^u} \\ 
			&\qquad\quad\, +\|\nabla w_h\|_{L^p}^2\|w_h\|_{L^q}\|v_h\|_{L^r}\big),
	\end{split}
\end{align}
where $1\leq p,\,q,\,r,\,s,\,t,\,u\leq \infty$ are another set of arbitrary numbers satisfying that
\begin{equation*}
	\tfrac2p+\tfrac1q+\tfrac1r = 1,\qquad \tfrac1s+\tfrac2t+\tfrac1u = 1.
\end{equation*}
\end{lemma}

\begin{proof}
	We begin with the proof of \cref{eq:lumperr1}. Henceforth, we denote by $I_h \coloneqq \sum_{K\in \mathcal{T}_h} I_K$ the $\mathcal T_h$-piecewise nodal interpolation, where $I_K$ is the local nodal interpolation on the element $K$. This operator is well-defined for $\mathcal T_h$-piecewise $H^2$-regular functions such as $V$. Using the triangle inequality, we obtain that
	\begin{align}
		\label{eq:trii}
		\begin{split}
			&|\ell(Vv_h,w_h)-(Vv_h,w_h)_{L^2}| \\
			&\qquad \leq  |\ell(I_hVv_h,w_h)-(I_hVv_h,w_h)_{L^2}| + |(I_hVv_h,w_h)_{L^2}-(Vv_h,w_h)_{L^2}|.
		\end{split}
	\end{align}
Classical approximation results for $I_h$ yield for the second term	that
 \begin{align*}
		|(Vv_h,w_h)_{L^2}-(I_hVv_h,w_h)_{L^2}| 
		&\lesssim h^2\|v_h\|_{L^v}\|w_h\|_{L^w}.
	\end{align*}

To estimate the first term on the right-hand side of~\cref{eq:trii}, we will use for each element a transformation to the reference simplex $\hat K$. For the simplex $K \in \mathcal T_h$, this transformation is an affine linear mapping given by $F_K \colon \hat K \to K,$ $x \mapsto B_Kx +b_K$, where $B_K \in \mathbb R^{d\times d}$ and $b_K \in \mathbb R^d$. We introduce for any simplex $K$ the functional
	\begin{align*}
		E_{K}\colon H^2(K)\to  \mathbb R,\quad v \to (v,1)_{L^2(K)} - \ell_K(v,1),
	\end{align*}
	which measures the mass lumping error. Note that, due to the continuous embedding $H^2(K)\hookrightarrow \mathcal C^0(\overline{K})$ this functional is well-defined and continuous. To estimate the norm of the functional 	$E_{K}$ we will perform a transformation to the reference simplex using the map $F_K$ and apply the Bramble-Hilbert lemma (see, e.g.,~\cite[Thm.~3.27]{knabner2021numerical}) to estimate the resulting functional $E_{\hat K}$ on the reference simplex. Note that in this proof we denote quantities that have been transformed to the reference simplex using a hat.  The change of variables formula for integrals then yields that
	\begin{align}
		\label{eq:transformation}
		E_K(v) = |\det(B_K)| E_{\hat K}(\hat v).
	\end{align}
To estimate $ E_{\hat K}(\hat v)$ we now apply the Bramble-Hilbert lemma. Since for all functions $\hat v \in \mathcal P^1(\hat K)$ it holds that $E_{\hat K}(\hat v) = 0$, we obtain the estimate
\begin{equation}
	\label{eq:bh}
	|E_{\hat K}(\hat v)| \leq C_\mathrm{BH} \|E_{\hat K}\||\hat v|_{H^2(\hat K)},
\end{equation}
where $C_\mathrm{BH}>0$ and $\|E_{\hat K}\|$ denotes the finite operator norm of the functional $E_{\hat K}$. 

Applying \cref{eq:transformation,eq:bh} for the particular  function~$v = q_h v_h w_h$, where we abbreviate $q_h\coloneqq I_hV$, yields that
\begin{align}
	\label{eq:tbc}
	E_K(q_h v_h w_h) 
	\leq C_\mathrm{BH}|\det(B_K)| \|E_{\hat K}\||\hat q_h\hat v_h \hat w_h|_{H^2(\hat K)}.
\end{align}
By the equivalence of norms in finite dimensions, we obtain that $|\hat q_h\hat v_h \hat w_h|_{H^2(\hat K)} \approx |\hat q_h\hat v_h \hat w_h|_{W^{2,1}(\hat K)}$. This result can then be used to continue \cref{eq:tbc} as follows
\begin{align*}
	&E_K(q_hv_hw_h)\\[-1.25ex]
	&\quad \lesssim  |\det(B_K)|\,|\hat q_h\hat v_h \hat w_h|_{W^{2,1}(\hat K)} = |\det(B_K)| \sum_{i,j=1}^{d}\|\hat \partial_j\hat \partial_i (\hat q_h\hat v_h \hat w_h)\|_{L^1(\hat K)}\\[-.5ex]
	&\quad \lesssim |\det(B_K)| \Big(\|\hat  q_h\|_{L^\infty(\hat K)}\|\hat \nabla \hat v_h\|_{L^p(\hat K)}\|\hat \nabla \hat w_h\|_{L^q(\hat K)} \\ &\qquad\quad + \|\hat \nabla \hat q_h\|_{L^6(\hat K)}\big(\|\hat \nabla v_h\|_{L^r(\hat K)}\|\hat w_h\|_{L^s(\hat K)}+\|\hat  v_h\|_{L^t(\hat K)}\|\hat \nabla \hat w_h\|_{L^u(\hat K)}\big)\Big)\\[1ex]
	&\quad \lesssim \|B_K\|^2\Big(\| q_h\|_{L^\infty(K)}\| \nabla v_h\|_{L^p( K)}\| w_h\|_{L^q( K)} \\ &\qquad\quad + \| \nabla q_h\|_{L^6( K)}\big(\| \nabla v_h\|_{L^r( K)}\| w_h\|_{L^s( K)}+\| v_h\|_{L^t(K)}\| \nabla w_h\|_{L^u( K)}\big)\Big),
\end{align*}
where we have used Hölder's inequality as well as \cite[Thm.~3.26]{knabner2021numerical} to transform the $L^p$-norms on the reference element back to the physical element. Further, we have used the local estimates
\begin{equation*}
	\|I_K V\|_{L^\infty(K)} \lesssim \|V\|_{H^2(K)},\qquad \|\nabla I_K V\|_{L^6(K)} \lesssim \|I_K V\|_{H^2(K)} \lesssim \|V\|_{H^2(K)}
\end{equation*}
which hold for all $K \in \mathcal T_h$. The estimate for the first term on the right-hand side of~\cref{eq:trii} can be concluded using the bound $\|B_K\| \lesssim h_K$ from \cite[Thm.~3.27]{knabner2021numerical} for the norm of the matrix $B_K$, and after summing over all elements $K \in \mathcal T_h$. 
Assertion~\cref{eq:lumperr1} then follows directly. The proof of \cref{eq:lumperr2} uses similar arguments and will be omitted for the sake of brevity. 
\end{proof}

The next lemma proves $L^\infty$-bounds for the  ground state approximations~$u_h$ defined in~\cref{eq:gsdisc}, and the ground state approximations obtained by the standard linear~FEM, denoted by $\hat u_h$.
  
\begin{lemma}[Uniform $L^\infty$-bounds]\label{lem:linftybounds}
Suppose that the energies $\mathcal E_h(u_h)$ and $\mathcal E(\hat u_h)$ are uniformly bounded. Then, the  $L^\infty$-bounds
	\begin{align*}
		\|u_h\|_{L^\infty} \lesssim 1,\qquad \|\hat u_h\|_{L^\infty} \lesssim 1
	\end{align*}
 hold uniformly in $h>0$.
\end{lemma}
\begin{proof}
	Let us first prove the $L^\infty$-bound for the ground state approximations $\hat u_h$ of the standard linear FEM. We define 
$\hat u_h^c \in H^1_0(\Omega)$ as the weak solution to 
	\begin{align}
		\label{eq:problem4uhhc}
	-\Delta \hat u_h^c = -V \hat u_h -\kappa|\hat u_h|^2\hat u_h + \lambda_h\hat u_h \eqqcolon f_h.
\end{align}
It holds that $f_h\in L^2(\Omega)$ with $\|f_h\|_{L^2}\lesssim 1$. Classical elliptic regularity theory on convex domains can then be used to show that $\hat u_h^c \in H^2(\Omega)\cap H^1_0(\Omega)$ with the estimate $\|\hat u_h^c\|_{H^2} \lesssim \|f_h\|_{L^2} \lesssim 1$. Note that the discrete function $\hat u_h \in V_h^0$ is the Galerkin approximation to $\hat u_h^c$ and therefore satisfies the classical error estimate
\begin{align*}
	h^{-1}\|\hat u_h^c - \hat u_h\|_{L^2} + \|\nabla(\hat u_h^c - \hat u_h)\|_{L^2} \lesssim h |\hat u_h^c|_{H^2}.
\end{align*}
Denoting by $I_h \colon \mathcal C^0(\overline{\Omega})\to V_h$ the nodal interpolation, we can estimate the $L^\infty$-norm of $\hat u_h$ using the triangle inequality as
\begin{align*}
	\|\hat u_h\|_{L^\infty} \lesssim \|\hat u_h - I_h\hat u_h^c\|_{L^\infty} + \|I_h \hat u_h^c\|_{L^\infty} \eqqcolon \Xi_1 + \Xi_2.
\end{align*}
The summand $\Xi_1$ can be estimated using a classical comparison result for $L^p$-norms of discrete functions and the approximation properties of the nodal interpolation~as
\begin{align*}
	\Xi_1 &\lesssim h^{-d/2} \|\hat u_h-I_h\hat u_h^c\|_{L^2} \leq h^{-d/2}\big(\|\hat u_h-\hat u_h^c\|_{L^2} + \|\hat u_h^c - I_h \hat u_h^c\|\big)\\
	&\lesssim h^{2-d/2}|\hat u_h^c|_{H^2} \lesssim 1.
\end{align*}
For the term $\Xi_2$, we obtain using the continuous embedding $H^2(\Omega) \hookrightarrow \mathcal C^0(\overline{\Omega})$ that
\begin{align*}
	\Xi_2 \leq \|\hat u_h^c\|_{L^\infty} \lesssim \|\hat u_h^c\|_{H^2} \lesssim 1.
\end{align*}
Combining the above bounds proves the uniform $L^\infty$-bound for $\hat u_h$. 

For proving the uniform $L^\infty$-bound for $u_h$, we need to derive a problem similar to \cref{eq:problem4uhhc} with an $L^2$-right-hand side. Note that the functional $$F(v_h) \coloneqq \ell(-Vu_h-|u_h|^2u_h + \lambda_hu_h,v_h)$$
is in the dual space of $V_h\subset L^2(\Omega)$, which means that, by the Riesz representation theorem, there exists $g_h \in V_h$ such that $	(g_h,v_h)_{L^2} = F(v_h)$ holds for all $v_h \in V_h$. It also yields the following bound for the $L^2$-norm of $g_h$:
\begin{align}
	\label{eq:tau}
	\|g_h\|_{L^2} = \sup_{v_h \in V_h^0} \frac{F(v_h)}{\|v_h\|_{L^2}}\lesssim \|-Vu_h - |u_h|^2u_h + \lambda_hu_h\|_\ell.
\end{align}
Here  we used that on the space $V_h$, the norm $\|\cdot\|_\ell$ is uniformly equivalent to the $L^2$-norm. Below, we will estimate the terms on the right-hand side of~\cref{eq:tau} separately. Using the assumed boundedness of the energies, the embedding $H^1(\Omega)\hookrightarrow L^q(\Omega)$ for $1 \leq q \leq 6$, and the norm equivalence of $\|\cdot\|_\ell$ and the $L^2$-norm on~$V_h$, we obtain the estimates $\|Vu_h\|_\ell \lesssim \|V\|_{L^\infty}\|u_h\|_\ell \lesssim \|V\|_{L^\infty}\|u_h\|_{L^2} \lesssim 1$ and $\lambda_h\|u_h\|_\ell \lesssim \lambda_h \|u_h\|_{L^2} \lesssim 1$. It remains to show the uniform boundedness of the the second term on the right-hand side of \cref{eq:tau}. Note that, similar to \cite[Thm.~3.46]{knabner2021numerical}, one can show the uniform equivalence of the $L^6$-norm and the norm defined by
\begin{align*}
	\|v_h\| \coloneqq \bigg(\sum_{K\in \mathcal{T}_h}\frac{|K|}{d+1}\sum_{j = 1}^{d+1}v_h^6(p_{\tau_K(j)})\bigg)^{1/6}
\end{align*}
on the space $V_h$. This yields the bound
\begin{align*}
	\||u_h|^2u_h\|_\ell^2 = \|u_h\|^6 \lesssim \|u_h\|_{L^6}^6 \lesssim 1,
\end{align*}
which, together with the previous estimates proves that $\|g_h\|_{L^2} \lesssim 1$. 
The right-hand side $g_h$ now takes the place of $f_h$ in \cref{eq:problem4uhhc}, and proceeding similarly as above for $\hat u_h$  gives the desired uniform $L^\infty$-bound for $u_h$.
\end{proof}

\bibliographystyle{alpha}
\bibliography{bib}
\end{document}